\documentclass{amsart}
\usepackage{amsthm}
\usepackage{amsmath}
\usepackage{amssymb}
\usepackage{url}
\usepackage{booktabs}
\usepackage{enumerate}

\newcommand{\ZZ}{\mathbf{Z}}
\newcommand{\QQ}{\mathbf{Q}}
\newcommand{\FF}{\mathbf{F}}

\begin{document}

\title{Irregular primes to two billion}

\author{William Hart}
\address[William Hart]{Technische Universit\"at Kaiserslautern, Fachbereich Mathematik, Postfach 3049, 67653 Kaiserslautern, Germany}
\email{goodwillhart@googlemail.com}

\author{David Harvey}
\address[David Harvey]{School of Mathematics and Statistics, University of New South Wales, Sydney NSW 2052, Australia}
\email{d.harvey@unsw.edu.au}

\author{Wilson Ong}
\address[Wilson Ong]{University of Cambridge, Department of Engineering, Information Engineering Division, Trumpington Street, Cambridge, CB2 1PZ, United Kingdom}
\email{wo218@cam.ac.uk}

%\date{}

% notational changes relative to the code:
%   gamma = g
%   omega = g_1
%   theta = g_2
%   gamma^(2mk) = h_1
%   gamma^(n i_2) = (h_2)^(i_2)
%   a_i = x_i

\maketitle

\begin{abstract}
We compute all irregular primes less than $2^{31} = 2\,147\,483\,648$.
We verify the Kummer--Vandiver conjecture for each of these primes, and we check that the $p$-part of the class group of $\QQ(\zeta_p)$ has the simplest possible structure consistent with the index of irregularity of $p$.
Our method for computing the irregular indices saves a constant factor in time relative to previous methods,
by adapting Rader's algorithm for evaluating discrete Fourier transforms.
\end{abstract}

\section{Introduction and summary of results}

For each of the $105\,097\,564$ odd primes less than $2^{31} = 2\,147\,483\,648$, we performed the following tasks:
\begin{enumerate}[(1)]
\item We computed the \emph{irregular indices} for $p$, that is, the integers $r \in \{2, 4, \ldots, p-3\}$ for which $B_r = 0 \pmod p$,
where $B_r$ is the $r$-th Bernoulli number.
A pair $(p, r)$, with $r$ as above, is called an \emph{irregular pair}, and such an integer $r$ is called an \emph{irregular index} for $p$.
The number of such $r$ is called the \emph{index of irregularity} of $p$, denoted $i_p$.
A prime $p$ is called \emph{regular} if $i_p = 0$, and \emph{irregular} if $i_p > 0$.
\item We verified the Kummer--Vandiver conjecture for $p$, which asserts that the class number of the maximal real subfield of $\QQ(\zeta_p)$ is not divisible by $p$, where $\zeta_p$ denotes a primitive $p$-th root of unity.
\item We verified that the Iwasawa invariants satisfy $\lambda_p = \nu_p = i_p$.
This implies that for all $n \geq 1$, the $p$-part of the class group of $\QQ(\zeta_{p^n})$ is isomorphic to $(\ZZ/p^n \ZZ)^{i_p}$.
\end{enumerate}
For further background on irregular primes, including their role in algebraic number theory and the historical connection to Fermat's Last Theorem, we refer the reader to \cite{Was-cyclotomic} and \cite{Lan-cyclotomic-combined}.

The tabulation of irregular pairs has been repeated by many people since the mid-19th century.
As shown in Table \ref{tab:timeline}, the search bound has increased dramatically, as more powerful computing hardware and increasingly sophisticated algorithms have been deployed.
The total running time of our computation was approximately 8.6 million core-hours (almost 1000 core-years).

\newcommand{\notreported}{\textsuperscript{\textdagger} }

\begin{table}[h]
\caption{Timeline of irregular prime search bounds. The symbol \textdagger{} indicates that a computation was reported but the corresponding list of irregular pairs was not published.}
\label{tab:timeline}
\begin{tabular}{rrl}
\toprule
Year       &           Bound   & Authors \\
\midrule
      1850 &       $p \leq 43$ & Kummer \cite{Kum-allgemeiner} \\
      1851 &                97 & Kummer \cite{Kum-memoire} \\
      1874 &               163 & Kummer \cite{Kum-1874} \\
      1930 &               199 & Stafford \& Vandiver \cite{SV-irregular} \\ % reported
      1937 &               613 & Vandiver \cite{Van-FLT} (many papers from 1930 to 1937) \\ % reported
      1954 &            2\,000 & Lehmer, Lehmer \& Vandiver \cite{LLV-FLT}       \\ % reported
      1954 &            2\,503 & Vandiver \cite{Van-attack} \\ % reported
      1955 &            4\,001 & Selfridge, Nicol \& Vandiver \cite{SNV-FLT} \\ % reported
      1963 &           10\,000 & Lehmer \notreported \cite{Leh-automation} \\
      1964 &           25\,000 & Selfridge \& Pollack \notreported (see \cite{Joh-irregular}) \\
      1970 &            5\,500 & Kobelev \cite{Kob-fermat} \\ % haven't managed to find a copy, but seems pretty clear that the list would have appeared in the paper
      1973 &            8\,000 & Johnson \cite{Joh-vanishing} \\ % reported
      1975 &           30\,000 & Johnson \cite{Joh-irregular} \\  % reported
      1976 &           32\,768 & Wada \notreported (see \cite{Wag-irregular}) \\
      1978 &          125\,000 & Wagstaff \cite{Wag-irregular} \\   % reported
      1987 &          150\,000 & Tanner \& Wagstaff \notreported \cite{TW-bernoulli}  \\
      1991 &          400\,000 & Sompolski \notreported \cite{Som-fermat} \\
      1992 &       1\,000\,000 & Buhler, Crandall \& Sompolski \notreported \cite{BCS-million}   \\
      1993 &       4\,000\,000 & Buhler, Crandall, Ernvall \& Mets\"ankyl\"a \notreported \cite{BCEM-four}     \\
      1996 &       8\,000\,000 & Shokrollahi \notreported \cite{Sho-eight}     \\
      2001 &      12\,000\,000 & Buhler, Crandall, Ernvall, Mets\"ankyl\"a \& Shokrollahi \cite{BCEMS-twelve}  \\
      2011 &     163\,577\,856 & Buhler \& Harvey \cite{BH-irregular} \\  % reported
      2016 &  2\,147\,483\,648 & this paper \\
\bottomrule
\end{tabular}
\end{table}

The complete table of irregular pairs for $p < 2^{31}$ is available at the second author's web page.
Previously, the largest known value of $i_p$ was $7$.
We found many new primes with $i_p = 7$, four primes with $i_p = 8$, namely
 \[ p = 381\,348\,997,\; 717\,636\,389,\; 778\,090\,129,\; 1\,496\,216\,791, \]
and exactly one prime with $i_p = 9$, namely $p = 1\,767\,218\,027$.
For this last $p$, we found that $B_r = 0 \pmod p$ for the following nine values of $r$:
\begin{multline*}
63\,562\,190,\; 274\,233\,542,\; 290\,632\,386,\; 619\,227\,758,\; 902\,737\,892, \\ 1\,279\,901\,568,\; 1\,337\,429\,618,\; 1\,603\,159\,110,\; 1\,692\,877\,044.
\end{multline*}

A standard heuristic commonly attributed to Lehmer and Siegel asserts that the quantity $i_p$ should follow a Poisson distribution with parameter $1/2$.
That is, the asymptotic density of primes for which $i_p = m$ should be equal to $P_m = e^{-1/2} / 2^m m!$.
Table \ref{tab:density} gives the exact counts for the primes up to $2^{31}$.
Here $N = 105\,097\,564$ is the total number of primes, and $N_m$ is the number of primes $p$ for which $i_p = m$.

\begin{table}[h]
\caption{Irregular index statistics for $p < 2^{31}$}
\begin{tabular}{crll}
\toprule
$m$ &   $N_m$       & $N_m/N$         & $P_m$           \\
\midrule
0   & 63\,751\,120  & 0.606590        & 0.606531        \\
1   & 31\,873\,681  & 0.303277        & 0.303265        \\
2   &  7\,963\,496  & 0.075772        & 0.075816        \\
3   &  1\,326\,171  & 0.0126185       & 0.0126361       \\
4   &     165\,211  & 0.00157198      & 0.00157951      \\
5   &      16\,410  & 0.000156141     & 0.000157951     \\
6   &       1\,384  & 0.000013169     & 0.000013163     \\
7   &           86  & 0.000000818     & 0.000000940     \\
8   &            4  & 0.0000000381    & 0.0000000588    \\
9   &            1  & 0.0000000095    & 0.0000000033    \\
\bottomrule
\end{tabular}
\label{tab:density}
\end{table}

At a glance, the data appears to fit the Poisson hypothesis very tightly.
We tested this formally by performing a chi-square goodness of fit test, combining the categories for $m \geq 7$ into a single bucket.
In other words, we computed the test statistic $X = \sum_{m=0}^7 (O_m - E_m)^2 / E_m$, where $O_m = N_m$ and $E_m = P_m N$ for $m = 0, \ldots, 6$, and $O_7 = N_7 + N_8 + N_9$, $E_7 = 1 - \sum_{m=0}^6 E_m$.
Under the Poisson hypothesis, and under suitable asymptotic normality assumptions, $X$ should follow a $\chi^2$ distribution with 7 degrees of freedom.
The observed value is $X = 13.807$, and the corresponding $p$-value under the $\chi^2$ distribution is $P(X \geq 13.807) = 0.0547$.
While this technically does not allow us to reject the Poisson hypothesis at the conventional 5\% significance level, the value is sufficiently extreme as to warrant further investigation.

We therefore repeated the analysis for the primes in each of the 16 subintervals $(2^{27}j, 2^{27}(j+1))$ for $j = 0, 1, \ldots, 15$.
For each subinterval, we counted the number of primes in that interval for which $i_p = m$, combining the values for $m \geq 6$ into a single bucket.
We then define a test statistic $X$ analogously to the previous paragraph, which should follow a $\chi^2$ distribution with 6 degrees of freedom.
The observed values $X_j$ are shown in Table \ref{tab:chisqr}, together with the corresponding $p$-values.
This table gives us much more confidence in the Poisson hypothesis; only the $j = 10$ interval has a relatively extreme $p$-value, and it is unsurprising for this to occur for one of the 16 intervals.
Despite this conclusion, we point out that it is still unknown whether there are infinitely many regular primes!

\begin{table}[h]
\caption{$\chi^2$ statistics for subintervals $(2^{27}j, 2^{27}(j+1))$}
\begin{tabular}{ccccccccc}
\toprule
$j$         & 0      & 1      & 2      & 3      & 4      & 5      & 6      & 7      \\
$X_j$       &  4.654 & 10.559 &  6.472 &  3.209 &  3.099 &  8.530 &  4.025 &  5.819 \\
$p$-value   &  0.589 &  0.103 &  0.372 &  0.782 &  0.796 &  0.202 &  0.673 &  0.444 \\
\midrule
$j$         & 8      & 9      & 10     & 11     & 12     & 13     & 14     & 15     \\
$X_j$       &  1.908 &  4.657 & 12.850 &  5.972 &  4.538 &  1.612 &  6.765 &  4.716 \\
$p$-value   &  0.927 &  0.588 &  0.045 &  0.426 &  0.604 &  0.952 &  0.343 &  0.580 \\
\bottomrule
\end{tabular}
\label{tab:chisqr}
\end{table}

The remainder of the paper is structured as follows.
In Sections \ref{sec:irregular} to \ref{sec:general} we explain a new method for determining the irregular indices for each $p$, based on Rader's algorithm for discrete Fourier transforms \cite{Rad-prime}.
This saves a constant factor in time over previous approaches, and was by far the most expensive part of our computation, consuming about 7.1 million core-hours.

In Section \ref{sec:checks} we briefly discuss the verification of the Kummer--Vandiver conjecture and the cyclotomic invariant checks.
These cost respectively 720{,}000 core-hours and 820{,}000 core-hours.
Washington has given heuristics suggesting that the number of counterexamples to the Kummer--Vandiver conjecture for $p < x$ should grow like $\log \log x$ \cite[\S8.3]{Was-cyclotomic}, and that the number of primes $p < x$ for which $\lambda_p \neq i_p$ should also grow like $\log \log x$ \cite[Appendix to Ch.~10]{Lan-cyclotomic-combined}.
Of course, $\log \log x$ increases very slowly, and these predictions are entirely consistent with our failure to observe a single exceptional event so far.

In Section \ref{sec:hardware} we describe the hardware that we used, and in Section \ref{sec:correctness} we discuss issues relating to the correctness of the results. Finally, in Section \ref{sec:history} we mention several interesting facts that we came across during preparation of this article, concerning the history of the computation of irregular primes, which do not appear to be well known.

\section{Computing the irregular indices --- basic framework}
\label{sec:irregular}

Let $p$ be an odd prime.
Our point of departure is the following well-known congruence.
Fix an integer $c$ in the interval $0 < c < p$, and for $0 < x < p$ define
\[ f_c(x) = \left\lfloor \frac{c \cdot (x/c \bmod p)}p \right\rfloor - \frac{c-1}2 \in \QQ, \]
where $x/c \bmod p$ means the smallest nonnegative integer congruent to $x/c$ modulo~$p$.
Then for even $r$ in the range $2 \leq r \leq p-3$ we have
\begin{equation}
\label{eq:voronoi}
\frac{c^r - 1}{r} B_{r} = \sum_{x=1}^{p-1} x^{r-1} f_c(x) \pmod p.
\end{equation}
This can be proved by reading the statement of Theorem 2.3 of \cite[Ch.~2]{Lan-cyclotomic-combined} modulo $p$ (as is done in the proof of Corollary 2 to that theorem).
Note that this congruence determines $B_r \pmod p$ if $c^r \neq 1 \pmod p$, but contains no useful information if $c^r = 1 \pmod p$.

We rewrite the congruence as follows.
First scale and reindex the Bernoulli numbers by setting $j = r-1$ and
\[ b_j = \frac{(c^{j+1} - 1)}{j+1} B_{j+1}. \]
Note that $j$ ranges over $1, 3, 5, \ldots, p-4$.
Let $\gamma$ be a generator of $(\ZZ/p\ZZ)^\times$, and reorder the sum by putting
 \[ a_i = f_c(\gamma^i), \qquad 0 \leq i \leq p-2. \]
Then the congruence becomes simply
\[ b_j = \sum_{i=0}^{p-2} \gamma^{ij} a_i \pmod p. \]
This expresses $b_1, b_3, \ldots, b_{p-4}$ as the odd-index components of the DFT (discrete Fourier transform)
of the vector $(a_0, a_1, \ldots, a_{p-2})$ over $\ZZ/p\ZZ$ with respect to the $(p-1)$-th root of unity $\gamma$.
For later convenience we will also define $b_{p-2}$ according to this formula (even though $B_{p-1}$ is not $p$-integral).

To take advantage of small factors of $p-1$, we decompose this DFT into shorter DFTs using the Good--Thomas method \cite{Goo-fourier,Tho-physics}.
Let $m$ be the $7$-smooth part of $(p-1)/2$, that is, the product of all powers of $2$, $3$, $5$ and $7$ appearing in $(p-1)/2$, and let $n = (p-1)/2m$.
In particular, $n$ is not divisible by $2$, $3$, $5$ or $7$, and is relatively prime to $2m$.
Thus we may reindex $i$ and $j$ by
\[ i = 2m i_1 + n i_2, \qquad j = 2m j_1 + n j_2, \]
where
\[ 0 \leq i_1, j_1 < n, \qquad 0 \leq i_2, j_2 < 2m, \qquad \text{$j_2$ odd}. \]
Put
\[ \omega = \gamma^{4m^2}, \qquad \theta = \gamma^{n^2}. \]
Note that $\omega$ has order $n$ in $(\ZZ/p\ZZ)^\times$, and $\theta$ has order $2m$.
Writing $b_{j_2,j_1} = b_j$ and $a_{i_2,i_1} = a_i$, we arrive at the two-dimensional DFT
\[ b_{j_2,j_1} = \sum_{i_1=0}^{n-1} \sum_{i_2=0}^{2m-1} \omega^{i_1 j_1} \theta^{i_2 j_2} a_{i_2,i_1}. \]

We may cut the number of terms in half by observing that $f_c(p-x) = -f_c(x)$ for all $x$.
This implies that $a_{i_2 + m, i_1} = -a_{i_2, i_1}$, and therefore that
\[ b_{j_2,j_1} = 2 \sum_{i_1=0}^{n-1} \sum_{i_2=0}^{m-1} \omega^{i_1 j_1} \theta^{i_2 j_2} a_{i_2,i_1}. \]
Our overall strategy will be to decompose this transform into $m$ \emph{horizontal} DFTs of length~$n$ over $\ZZ/p\ZZ$, given by
 \[ d_{i_2,j_1} = 2 \sum_{i_1=0}^{n-1} \omega^{i_1 j_1} a_{i_2,i_1} \]
for $0 \leq i_2 < m$, followed by $n$ \emph{vertical} DFTs of length $m$ over $\ZZ/p\ZZ$, namely
 \[ b_{j_2,j_1} = \sum_{i_2=0}^{m-1} \theta^{i_2 j_2} d_{i_2,j_1} \]
for $0 \leq j_1 < n$.

With this high-level description of the algorithm in place, we may already deduce a lower bound for the overall memory requirements.
Indeed, after the horizontal DFTs, the entire intermediate array $d_{i_2,j_1}$ is stored in memory.
Assuming 32 bits per coefficient (since $p < 2^{31}$), for $p$ near the top of our search range this corresponds to about 4 GB.
In fact, as we will see below, for many $p$ this is a serious underestimate;
the peak memory usage, which occurs during the horizontal DFT stage, can be much higher.

To conclude this section, we briefly describe how we perform the vertical DFTs.
Each one has $m$ input coefficients, but $\theta$ has order $2m$, and we only want the odd-index components of the DFT.
We can convert this to a plain vanilla DFT by writing $j_2 = 2j' + 1$; we then want to evaluate
 \[ \sum_{i_2=0}^{m-1} (\theta^2)^{i_2 j'} (\theta^{i_2} d_{i_2,j_1}) \]
for $0 \leq j' < m$.
Given the inputs $d_{i_2, j_1}$ (i.e., the outputs of the horizontal DFT stage),
we simply twist them by the coefficients $\theta^{i_2}$,
and then perform an FFT of length $m$ with respect to the $m$-th root of unity $\theta^2$,
working directly over $\ZZ/p\ZZ$.
Since $m$ was constructed to be $7$-smooth, it is straightforward to use the Cooley--Tukey FFT algorithm to decompose the transform into small hard-coded transforms of length $2$, $3$, $5$ and $7$.
This is all done in-place, so the memory required in this phase is at most 4 GB.

For most primes $p < 2^{31}$, the factor $m$ is quite small: the median value of $m$ is 8, and for 95\% of primes it is less than $500$.
Consequently, our code for the vertical DFTs is fairly simple: we focused our optimisation efforts on fast modular arithmetic,
using standard methods to reduce the number of arithmetic operations,
and we did not worry too much about locality, which becomes much more important for the horizontal DFTs.
Our implementation is certainly suboptimal for those rare primes~$p$ for which $m$ is very large (although it is still reasonably fast).

\section{Horizontal DFTs --- the umbrella algorithm}

Continuing with the notation established above, let us fix $i_2$ and consider the horizontal DFT for the $i_2$-th row.
Setting
 \[ d_\ell = d_{i_2,\ell}, \qquad a_k = a_{i_2,k} = f_c(\gamma^{n i_2} (\gamma^{2m})^k), \]
the DFT in question may be written as
 \[ d_\ell = 2 \sum_{k=0}^{n-1} \omega^{k\ell} a_k, \qquad 0 \leq \ell < n. \]

We used several methods to compute this DFT.
The method chosen for each $p$ depends heavily on the factorisation of $p-1$ and on the number-theoretic properties of the ring $\ZZ/p\ZZ$.
In subsequent sections we present fast methods that work for most primes; a prime for which they do not work is said to be \emph{rejected}.
The rejected primes are handled by the \emph{umbrella algorithm}, which we explain in this section.
It is slower, but works for any prime.
It is essentially equivalent to the ``Voronoi identity method'' of \cite{BH-irregular}
(and is closely related to the ``root finding method'' of \cite{BCEMS-twelve}).

Take $c = \gamma$, so that \eqref{eq:voronoi} determines $B_r \pmod p$ for all $r$.
We will apply Bluestein's method \cite{Blu-dft} to convert the horizontal DFT to a convolution problem.
Put $\xi = \gamma^{2m^2}$ (so that $\xi^2 = \omega$); then
 \[ d_\ell = 2 \xi^{\ell^2} \sum_{k=0}^{n-1} \xi^{-(\ell-k)^2} \xi^{k^2} a_k. \]
One checks that $\xi^{-(k+n)^2} = \xi^{-k^2}$ for all $k$.
Thus, if we define two polynomials
 \[ U(X) = \sum_{k=0}^{n-1} \xi^{k^2} a_k X^k \in \FF_p[X], \qquad V(X) = \sum_{k=0}^{n-1} \xi^{-k^2} X^k \in \FF_p[X], \]
and put
 \[ W(X) = U(X) V(X) = \sum_{\ell=0}^{2n-1} W_\ell X^\ell, \]
then we find that $d_\ell = 2 \xi^{\ell^2} (W_\ell + W_{\ell + n})$.
Consequently the problem boils down to computing the product $U(X) V(X)$ of two polynomials of degree less than $n$ in $\FF_p[X]$,
plus various pre- and post-processing operations.

We briefly discuss these auxiliary operations first.
These amount to computing the sequences $\xi^{\ell^2}$, $\xi^{-\ell^2}$, $\gamma^{n i_2} (\gamma^{2m})^k$ and $a_k$, plus $O(n)$ multiplications and additions in $\FF_p$.
This can all be achieved in no more than $O(n)$ word operations.

Now we consider the polynomial multiplication step.
If we are lucky enough that $p-1$ has only small prime factors, then this could be done efficiently via FFTs working directly over $\ZZ/p\ZZ$.
However, this approach is viable only for a very thin set of primes, so we used instead the following method, which is applicable to any~$p$.

We first lift the multiplication problem to $\ZZ[X]$.
Let $\tilde U(X)$ and $\tilde V(X)$ denote lifts of $U(X)$ and $V(X)$ to $\ZZ[X]$,
with coefficients taken in the interval $-p/2 < x < p/2$, and let $\tilde W(X) = \tilde U(X) \tilde V(X) \in \ZZ[X]$.
The degrees of $\tilde U(X)$ and $\tilde V(X)$ are at most $n - 1 < p/2$, so the coefficients of
$\tilde W(X)$ are bounded in absolute value by $p^3/8$.
For $p$ near the top of our search range, namely $p \sim 2^{31}$, we find that the coefficients of $\tilde W(X)$ have at most $91$ bits (i.e., 90 bits plus a sign bit).

As 91 bits exceeds the 64-bit word size of the machines available to us, we employ a multimodular approach.
We choose two 62-bit primes $q_1$ and $q_2$, with $q_i - 1$ divisible by a large power of two.
These primes are fixed at the beginning, independently of $p$.
We compute the product $\tilde U(X) \tilde V(X)$ first modulo $q_1$ and then modulo $q_2$, using FFTs over $\ZZ/q_i \ZZ$.
We then combine the results using the Chinese remainder theorem to deduce $\tilde W(X)$,
and finally reduce the result modulo $p$ to deduce $W(X)$.
Two primes suffice as $62 + 62 = 124 > 91$.

(Actually, because of cancellation between positive and negative terms, the coefficients of $\tilde W(X)$ are usually closer to 75 bits in size. Unfortunately, we cannot give a provable bound.
It does not help anyway, because $75$ bits is still too large for single-word arithmetic.)

The overall running time is dominated by the cost of these FFTs, and we spent considerable effort optimising them.
Our optimisation techniques may be summarised as follows.
For the low-level modular arithmetic, we use the method of \cite{Har-ntt} to reduce the number of modular reductions performed; this explains our choice of 62-bit primes rather than primes of 63 or 64 bits.
We use truncated Fourier transforms \cite{vdH-TFT-apps} to avoid power-of-two jumps in the running time; that is, taking $N$ to be a suitably large power of two, instead of evaluating at all $N$-th roots of unity, we evaluate only on a subset of those roots large enough to determine the polynomial product of interest.
We aggressively use array decompositions of FFTs, adapted to the truncated case \cite{Har-cachetft}, to improve locality.
Finally, we use OpenMP throughout, including within the FFTs themselves, to take advantage of multiple processor cores.

Note that $V(X)$ does not depend on the row index $i_2$.
To exploit this, we first compute $V(X)$ and then the transforms of $\tilde V(X)$ modulo $q_1$ and $q_2$.
Then for each row, we compute $U(X)$, then the transforms of $\tilde U(X)$, then multiply these by the transforms of $\tilde V(X)$ and
perform the inverse transforms to deduce $\tilde W(X)$.
Altogether we perform $4m + 2$ transforms of length approximately~$2n$.
The worst case is $m = 1$, in which we perform six transforms of length roughly~$p$, and the peak memory usage is roughly $30p$ bytes.
For $p \sim 2^{31}$ this amounts to 60~GB. This was the largest memory footprint encountered for any single prime $p$ considered in this paper.

\section{Horizontal DFTs --- Rader's method with one prime}
\label{sec:oneprime}

The main drawback of the umbrella algorithm of the previous section is the coefficient growth encountered when lifting the polynomial product to $\ZZ[X]$.
This has been a perennial problem in the history of irregular prime computations.
All known asymptotically fast algorithms for computing irregular indices ultimately reduce to the problem of multiplying polynomials in $\FF_p[X]$ of degree proportional to $p$.
After lifting to $\ZZ[X]$, this leads to a bound of the form $O(p^3)$ for the size of the coefficients of the product: two factors of $p$ arise from the size of the input coefficients, and one factor of $p$ from the length of the polynomials.
For example, in 1992, when Buhler, Crandall and Sompolski determined the irregular primes up to one million \cite{BCS-million},
they reported using two 30-bit primes (on 32-bit hardware) to recover coefficients of size roughly $(10^6)^3 \sim 2^{60}$.

We now describe a new approach, based on \emph{Rader's method} \cite{Rad-prime}, that reduces this bound from $O(p^3)$ to $O(p^2)$ for most primes.
In this section we focus on the easiest case in which $n$ is itself \emph{prime}; this restriction will be dropped in the following sections.

Recall the expression for the horizontal DFT for the $i_2$-th row:
 \[ d_\ell = 2 \sum_{k=0}^{n-1} \omega^{k\ell} a_k, \qquad 0 \leq \ell < n. \]
The case $\ell = 0$ will be handled separately; for this we have simply
 \[ d_0 = 2 \sum_{k=0}^{n-1} a_k. \]
Now let $z$ be a generator of $(\ZZ/n\ZZ)^\times$, and for $1 \leq \ell < n$ write $\ell = z^t \pmod n$, where $0 \leq t < n-1$.
Similarly for $1 \leq k < n$ write $k = z^{-s} \pmod n$, where $0 \leq s < n-1$.
The transform thereby becomes
 \[ d_{z^t} = 2 a_0 + 2 \sum_{s=0}^{n-2} \omega^{z^{t-s}} a_{z^{-s}}. \]
The last sum may be recognised as a cyclic convolution of length $n-1$.
Indeed, if we define
 \[ U(X) = \sum_{s=0}^{n-2} a_{z^{-s}} X^s \in \FF_p[X], \qquad V(X) = \sum_{s=0}^{n-2} \omega^{z^s} X^s \in \FF_p[X], \]
and
 \[ W(X) = U(X) V(X) = \sum_{t=0}^{2n-3} W_t X^t, \]
then we have
 \[ d_{z^t} = 2 a_0 + 2 (W_t + W_{t+n-1}). \]

We now proceed just as in the umbrella algorithm, by lifting the problem to multiplication in $\ZZ[X]$.
However, there is a crucial difference.
In the umbrella algorithm, the coefficients of $\tilde U(X)$ and $\tilde V(X)$ were random-looking integers in the interval $-p/2 < x < p/2$.
In the present situation, we may take the coefficients of $\tilde U(X)$ to be \emph{small} integers, provided that we choose $c$ appropriately.

Let us examine the cases $c = 2$ and $c = 3$.
For $c = 2$ we have
\[ f_2(x) = \begin{cases} -1/2 & x = 0 \pmod 2, \\ 1/2 & x = 1 \pmod 2, \end{cases} \]
so $2a_{z^{-s}} = \pm 1$ for each $s$.
Let $\tilde U(X)$ be a lift of $2 U(X)$ to $\ZZ[X]$ with coefficients in $\{-1, 1\}$, and let $\tilde V(X)$ be a lift of $V(X)$ with coefficients in the interval $-p/2 < x < p/2$.
The problem is now to compute the product $\tilde W(X) = \tilde U(X) \tilde V(X)$ in $\ZZ[X]$.
The coefficients of $\tilde W(X)$ are bounded in absolute value by $(p/2) n \leq p^2/4$.

Now consider the case $c = 3$. We have
 \[ f_3(x) = \begin{cases} -1 & x = 0 \pmod 3, \\ 1 & x = p \pmod 3, \\ 0 & x = -p \pmod 3. \end{cases} \]
We take $\tilde U(X)$ to be a lift of $U(X)$ with coefficients in $\{-1, 0, 1\}$, and take $\tilde V(X)$ and $\tilde W(X)$ as above.
Again the coefficients of $\tilde W(X)$ are bounded in absolute value by $p^2/4$.

In both cases the coefficients are bounded by $p^2/4 < 2^{60}$ for the range of primes under consideration. Thus to compute $\tilde W(X)$ it suffices to compute it modulo a single 62-bit prime $q$ (indeed, 61 bits would be enough). Consequently we save a factor of roughly two compared to the umbrella algorithm.

(Again, due to cancellation, it is likely that the coefficients of $\tilde W(X)$ are closer to 44 bits.
It is thus conceivable that one could use complex FFTs with double-precision floating point arithmetic, which are likely to be faster than FFTs over $\ZZ/q\ZZ$.
We decided in the end to take the modular path, mainly because it did not seem possible to obtain a provable bound, or to obtain sufficiently tight provable error bounds for the floating-point FFTs.)

We present a brief toy example to illustrate the advantage of the Rader algorithm over the umbrella algorithm.
Let $p = 131$, so that $m = 5$ and $n = 13$, and take $\gamma = 2$ and $c = 2$.
Consider the horizontal DFT corresponding to $i_2 = 3$.
In the umbrella algorithm, we need to compute the product of
\begin{align*}
 \tilde U(X) & = 40X^{12} + 16X^{11} - 43X^{10} - 26 X^9 + 9X^8 - 34X^7 \\
             & \qquad\qquad\qquad + 34X^6 - 9X^5 - 26 X^4 - 43X^3 + 16X^2 - 40X + 65, \\
 \tilde V(X) & = -18X^{12} + 45X^{11} - 32X^{10} + 63X^9 - 51X^8 + 52X^7 \\
             & \qquad\qquad\qquad + 52X^6 - 51X^5 + 63X^4 - 32X^3 + 45X^2 - 18X + 1.
\end{align*}
In the Rader version, we need to compute the product of
\begin{align*}
 \tilde U(X) & = -X^{11} - X^{10} - X^9 + X^8 - X^7 + X^6 - X^5 - X^4 + X^3 + X^2 + X - 1, \\
 \tilde V(X) & = -51 X^{11} + 39X^{10} + 63X^9 + 60X^8 + 45X^7 + 62X^6 \\
             & \qquad\qquad\qquad - 18X^5 - 47X^4 + 52X^3 - 24X^2 - 32X - 19.
\end{align*}

In our implementation we always choose $c = 2$ or $c = 3$.
One unpleasant byproduct of this choice is that $c$ may not generate $(\ZZ/p\ZZ)^\times$, so that \eqref{eq:voronoi} may not determine $B_r$ modulo $p$ for all $r$.
More precisely, writing $\alpha_c$ for the order of $c$ in $(\ZZ/p\ZZ)^\times$, equation \eqref{eq:voronoi} determines $B_r$ except for those $r$ divisible by $\alpha_c$.

We handle this as follows.
Given $p$ as input, we first try $c = 2$.
If $\alpha_2$ is sufficiently large, say $\alpha_2 > (p-1)/100$, we use the method described above, with $c = 2$, to compute $B_r$ for those $r$ not divisible by $\alpha_2$.
We then compute the missing values of $B_r$ separately, using a method described in the next paragraph.
If the order of~$2$ is too small, we try $c = 3$ instead.
If $\alpha_3 > (p-1)/100$, we proceed as before, with $c = 3$.
If we are unlucky, both $2$ and $3$ will have small order; in this case, $p$ is rejected.

Now we explain how to recover $B_r \pmod p$ for the missing values of $r$.
Put $\alpha = \alpha_c$.
There are $N - 1$ missing values, where $N = (p-1)/\alpha \leq 100$, namely $r = j\alpha$ for $j = 1, 2, \ldots, N-1$.
(Some of these $r$ may be odd, and could be skipped if desired.)
Applying \eqref{eq:voronoi} with $c = \gamma$, we obtain
\begin{align*}
 \frac{\gamma^{j\alpha} - 1}{j\alpha} B_{j\alpha}
   & = \sum_{i=0}^{p-2} (\gamma^\alpha)^{ji} \gamma^{-i} f_\gamma(\gamma^i) = \sum_{k=0}^{N-1} (\gamma^\alpha)^{jk} \sum_{i = k \bmod N} \gamma^{-i} f_\gamma(\gamma^i) \pmod p.
\end{align*}
This expresses the missing $B_{j\alpha}$ in terms of a DFT of length $N$ with respect to the $N$-th root of unity $\gamma^\alpha$.
The inputs to the DFT, i.e., the $N$ sums over $i$, are easily computed in $O(p)$ word operations.
As $N$ is very small, the cost of the DFT itself is insignificant compared to the main part of the algorithm.

It is certainly possible to throw other small values of $c$ into the mix, such as $c = 5$, and this would reduce the proportion of rejected primes.
The tradeoff is increased code complexity.
We found that using $c = 2$ and $c = 3$ strikes a reasonable balance.

Finally, we point out one additional complication that arises in applying Rader's method, namely, the problem of computing the polynomials $U(X)$ and $V(X)$ efficiently.
For $V(X) = \sum_{s=0}^{n-2} \omega^{z^s} X^s$, one simple approach is to notice that $\omega^{z^s} = (\omega^{z^{s-1}})^z$, so that each term may be computed from the previous one by a single $z$-th power modulo $p$.
This works well if $z$ is small, say $z = 2$ or $z = 3$.
Unfortunately, it may well happen that neither $2$ nor $3$ generates $(\ZZ/n\ZZ)^\times$.

As a workaround, we instead try to choose $z$ so that $z^M = y \pmod n$, where $y = 2$ or $y = 3$, and where $M$ is small, say $M \leq 100$.
Assuming this can be done, we compute $V(X)$ as follows: we first compute $\omega^{z^s}$ for $0 \leq s < M$ directly, and then repeatedly evaluate $\omega^{z^s} = (\omega^{z^{s-M}})^y$ for $s \geq M$.
This reduces the amount of work to roughly one ($y = 2$) or two ($y = 3$) multiplications modulo $p$ per coefficient of $V(X)$.

To select $y$, we use a strategy similar to that used earlier to choose $c$.
We first check the order of $2$ in $(\ZZ/n\ZZ)^\times$; if it is sufficiently large, we take $y = 2$. Otherwise we examine the order of $3$; if it is large, we take $y = 3$. Assuming one of these options works, it is then straightforward (since $M$ is small) to apply an analogue of the Tonelli--Shanks algorithm to find a suitable $z$ satisfying $z^M = y$. If neither option works, we reject this prime.

Analogous remarks apply to computing the coefficients of $U(X)$.
Note that the inner loop for computing $U(X)$ depends on both $c$ and $y$.
Our code contains four versions of this loop, one for each choice of $(c, y)$.

\section{Horizontal DFTs --- Rader's method with two primes}
\label{sec:twoprimes}

Let us suppose now that $n$ is a product of \emph{two distinct primes}, say $n = n_1 n_2$.
Consider again the horizontal DFT for the $i_2$-th row:
 \[ d_\ell = 2 \sum_{k=0}^{n-1} \omega^{k\ell} a_k, \qquad 0 \leq \ell < n. \]
There are four types of $\ell$ we must consider.
The first case is where $\ell = 0$; as before we simply have $d_0 = 2 \sum_{k=0}^{n-1} a_k$.
The second case is where $\ell$ is divisible by $n_2$ but not by $n_1$, say $\ell = n_2 \ell'$ where $1 \leq \ell' < n_1$.
Then the expression for $d_\ell$ becomes
\begin{equation}
\label{eq:n1-first}
 d_\ell = 2\sum_{k=0}^{n_1 - 1} (\omega^{n_2})^{k\ell'} a'_k,
\end{equation}
where
 \[ a'_k = \sum_{k'=0}^{n_2-1} a_{n_1 k' + k}, \qquad 0 \leq k < n_1. \]
This may be regarded as a DFT of length $n_1$ of the $a'_k$, with respect to the $n_1$-th root of unity $\omega^{n_2}$.
Since $n_1$ is prime, we may apply Rader's method to reduce this to a cyclic convolution of length $n_1 - 1$.
The input coefficients are larger than in the previous section; indeed, taking $c = 2$ or $c = 3$ as before, we find that $|a'_k| \leq n_2$ rather than $|a'_k| \leq 1$.
On the other hand, the convolution length is shorter.
These forces exactly counteract, so that we can still perform the multiplication in $\ZZ[X]$ using one 62-bit prime $q$.
The details are omitted.

The third case is where $\ell$ is divisible by $n_1$ but not $n_2$.
This is handled exactly as above, with the roles of $n_1$ and $n_2$ switched.

The fourth case is where $\ell$ is divisible by neither $n_1$ nor $n_2$, i.e., where $(\ell,n) = 1$.
Of course, for most primes $p$, almost all $\ell$ fall into this case.
We split the sum into four sums, according to the divisibility of $k$ by $n_1$ and $n_2$:
\begin{equation}
\label{eq:n1-second}
 d_\ell = 2a_0 + 2\sum_{k'=1}^{n_1 - 1} (\omega^{n_2})^{k' \ell} a_{n_2 k'} + 2 \sum_{k'=1}^{n_2 - 1} (\omega^{n_1})^{k' \ell} a_{n_1 k'} + 2 \sum_{\substack{1 \leq k < n \\ (k,n) = 1}} \omega^{k\ell} a_k.
\end{equation}
The second and third terms reduce to DFTs of length $n_1$ and $n_2$, and these are (yet again) converted to cyclic convolutions via Rader's method.

%Note that so far we have evaluated two convolutions of length $n_1 - 1$: one for the sum in \eqref{eq:n1-first}, and one for the second term in \eqref{eq:n1-second}.
%We have also evaluated two convolutions of length $n_2 - 1$, and performed $O(p)$ additional work, such as computing the $a'_k$.
%Typically $n_1$ and $n_2$ are much smaller than $n$, so the cost of these tasks is negligible.

We turn now to the fourth sum, which is where the bulk of the computation occurs.
The most natural generalisation of the algorithm of Section \ref{sec:oneprime} is a \emph{two-dimensional} variant of Rader's method, which runs as follows.
Consider the isomorphism $(\ZZ/n\ZZ)^\times \cong (\ZZ/n_1\ZZ)^\times \times (\ZZ/n_2\ZZ)^\times$.
Choose $z_i$ to be a generator of the $(\ZZ/n_i\ZZ)^\times$ part, i.e., $z_1$ is a generator modulo $n_1$, and $z_1 = 1 \pmod{n_2}$, and similarly for $z_2$.
Every $\ell \in (\ZZ/n\ZZ)^\times$ can be written uniquely in the form $\ell = z_1^{t_1} z_2^{t_2}$ for $0 \leq t_i < n_i - 1$.
Similarly, we write $k = z_1^{-s_1} z_2^{-s_2}$ for $0 \leq s_i < n_i - 1$.
Then the sum becomes
 \[ \sum_{s_1=0}^{n_1-2} \sum_{s_2=0}^{n_2-2} \omega^{z_1^{t_1 - s_1} z_2^{t_2 - s_2}} a_{z_1^{-s_1} z_2^{-s_2}}. \]
This is a two-dimensional cyclic convolution of size $(n_1 - 1) \times (n_2 - 1)$, i.e., evaluating it is equivalent to computing a certain product in $\FF_p[X_1,X_2]/(X_1^{n_1 - 1} - 1, X_2^{n_2 - 1} - 1)$.
This could be achieved by lifting to $\ZZ$, computing the product in $\ZZ[X_1, X_2]$ (which can be done modulo a suitable 62-bit prime $q$, just as in Section \ref{sec:oneprime}), and then performing the cyclic reductions.
Unfortunately, we encounter the following problem: the product in $\ZZ[X_1, X_2]$ has degree roughly $2n_1$ in $X_1$ and $2n_2$ in $X_2$, so the total transform size is about $4n$.
In other words, we have zero-padded in two dimensions instead of one, and this essentially doubles the running time relative to Section \ref{sec:oneprime}.
To remove this undesirable factor of two, we must work harder.

First we will construct factorisations
 \[ n_1 - 1 = d_1 e_1, \qquad n_2 - 1 = d_2 e_2, \]
where $d_1$ and $d_2$ are as large as possible, subject to the constraint that $(d_1, d_2) = 1$.
To do this, we consider in turn each prime divisor $\pi$ of $(n_1 - 1)(n_2 - 1)$.
Suppose that $\pi^{a_i}$ is the exact power of $\pi$ dividing $n_i - 1$.
If $a_1 \geq a_2$ we include $\pi^{a_1}$ in $d_1$ and $\pi^{a_2}$ in $e_2$; otherwise we place $\pi^{a_1}$ in $e_1$ and $\pi^{a_2}$ in $d_2$.
(For example, if $n_1 = 10459$ and $n_2 = 19249$, then $n_1 - 1 = 2 \cdot 3^2 \cdot 7 \cdot 83$ and $n_2 - 1 = 2^4 \cdot 3 \cdot 401$, so we put $d_1 = 3^2 \cdot 7 \cdot 83$, $e_1 = 2$, $d_2 = 2^4 \cdot 401$, $e_2 = 3$.)

The result is that $(\ZZ/n\ZZ)^\times \cong (\ZZ/d_1 d_2 \ZZ) \oplus (\ZZ/e_1\ZZ) \oplus (\ZZ/e_2\ZZ)$.
Now we can proceed just as before, using this alternative decomposition of $(\ZZ/n\ZZ)^\times$.
We select generators, say $u_0, u_1, u_2 \in (\ZZ/n\ZZ)^\times$, of the three cyclic subgroups.
The sum then takes the form
 \[ \sum_{s_0=0}^{d_1 d_2 - 1} \sum_{s_1=0}^{e_1 - 1} \sum _{s_2=0}^{e_2 - 1} \omega^{u_0^{t_0 - s_0} u_1^{t_1 - s_1} u_2^{t_2 - s_2}} a_{u_0^{-s_0} u_1^{-s_1} u_2^{-s_2}}. \]
This is a three-dimensional cyclic convolution of size $d_1 d_2 \times e_1 \times e_2$, and so it is equivalent to computing a product in $\FF_p[X_0, X_1, X_2]/(X_0^{d_1 d_2} - 1, X_1^{e_1} - 1, X_2^{e_2} - 1)$.
We perform this product by lifting to $\ZZ[X_0, X_1, X_2]/(X_1^{e_1} - 1, X_2^{e_2} - 1)$.
Note that we have \emph{retained} the cyclic structure with respect to $e_1$ and $e_2$.
Next we choose some 62-bit prime $q$ and perform the product in $\FF_q[X_0, X_1, X_2]/(X_1^{e_1} - 1, X_2^{e_2} - 1)$.
The crucial observation is that \emph{if $e_1$ and $e_2$ are sufficiently small}, then we can choose $q$ so that $q = 1 \pmod{e_1}$ and $q = 1 \pmod{e_2}$.
Then $\FF_q$ contains appropriate roots of unity to handle the cyclic products of length $e_1$ and $e_2$ directly, without zero-padding.
Of course we still require that $q = 1$ modulo a large power of two, to handle the FFT in the ``$d_1 d_2$ dimension'', and we perform the usual zero-padding in that dimension.

If $e_1$ and $e_2$ are not small (or smooth) enough, we reject this prime.
This occurs for example if $n_1 - 1$ and $n_2 - 1$ have a large prime factor in common.

All of the complications and optimisations mentioned for the single prime case (Section \ref{sec:oneprime}) still apply in this setting.
The main difference is notational complexity; we omit the details.

\section{Horizontal DFTs --- Rader's method in general}
\label{sec:general}

In principle, it is straightforward to generalise the algorithm of Section \ref{sec:twoprimes} to those $n$ having three or more prime factors.
For example, if $n = n_1 n_2 n_3$, then the problem reduces to a collection of ``small'' convolutions of sizes $n_i - 1$ and $(n_i - 1) \times (n_j - 1)$, together with one large four-dimensional convolution of size $(d_1 d_2 d_3) \times e_1 \times e_2 \times e_3$.

We implemented this for the cases where $n$ is a product of up to four distinct prime factors.
Beyond this, it does not seem to be worth the trouble, for a number of reasons.
First, there are simply not that many primes $p$ for which $n$ has five or more factors, so overall it does not hurt much to use the umbrella algorithm for these.
Second, as the number of factors of $n$ increases, we tend to reject a larger proportion of primes.
Third, the variants with more prime factors are simply less efficient, because of additional overheads in all stages of the algorithm.
Presumably, for a large enough search bound, it would become worthwhile to implement the five-prime case, and it would also be profitable to tweak the parameters so as to reduce the rejection rates.

Another case that we ignore is where $n$ has a \emph{repeated} prime factor $\geq 11$.
For example, this occurs for any prime $p = 1 \pmod{11^2}$.
The methods of Section \ref{sec:oneprime} and \ref{sec:twoprimes} could be modified to handle this case, but we did not bother because there did not seem to be enough primes in this category to justify the effort.
An alternative approach would be to redefine $m$ to be the $11$-smooth part of $(p-1)/2$, at the cost of pushing more work into the vertical DFT stage.

Table \ref{tab:numprimes} shows the number of odd primes $p < 2^{31}$ for which $n$ was a product of one, two, three or four distinct prime factors, and the number of primes fitting into none of those categories.
It also shows, for each category, the number of primes that were rejected (for any reason).

\begin{table}[h]
\caption{Classification of primes $p < 2^{31}$}
\begin{tabular}{lrrrr}
\toprule
Factorisation     & Initial number       &  Proportion of     & Number of       & Rejection         \\
of $n$            & of primes            &  all primes (\%)   & rejected primes & rate (\%)         \\
\midrule
$n_1$             &        26\,129\,901  &  24.9              &      8\,933     &   0.03            \\
$n_1 n_2$         &        46\,221\,996  &  44.0              & 1\,304\,998     &   2.8\phantom{0}  \\
$n_1 n_2 n_3$     &        24\,522\,583  &  23.3              & 1\,628\,092     &   6.6\phantom{0}  \\
$n_1 n_2 n_3 n_4$ &         4\,574\,382  &   4.4              &    854\,873     &  18.7\phantom{0}  \\
All other types   &         3\,648\,702  &   3.5              & 3\,648\,702     & 100.0\phantom{0}  \\
\midrule
Total             &       105\,097\,564  & 100.0              & 7\,481\,598     &   7.1\phantom{0}  \\
\bottomrule
\end{tabular}
\label{tab:numprimes}
\end{table}

Table \ref{tab:performance} shows the time and memory use of our implementation under various scenarios.
These tests were run on a single 2.6 GHz Intel Xeon E5-2650v2 processor on the \textsc{Katana} cluster (see Section \ref{sec:hardware}).
For each of the five algorithms, and for two choices of $m$, we selected a typical prime near the top of the search range, and we tested both the single-core performance and also the performance when running 8 cores simultaneously.
The time includes the cost of verifying the checksum discussed in Section \ref{sec:correctness}.
Regarding $m$, the worst possible value is $m = 1$, because there is no opportunity for FFT reuse; conversely, the table illustrates the benefit obtained in the case of the more typical value $m = 8$.

As the table shows, our code scales quite well up to 8 cores, at least on this machine; in all rows of the table the speedup is at least a factor of 7 compared to the single-core performance.
In reality, we did not often need to run the code on this many cores.
The majority of primes did not need too much RAM, and for these primes we preferred to use as few cores as possible; two or four cores was typical.
It was only for primes requiring large amounts of RAM, say between 30 GB and 60 GB, that we were forced to use larger core counts.

\begin{table}[h]
\caption{Performance of various algorithms for determining irregular indices. Time is wall clock time in seconds. Memory is peak RAM usage in GB.}
\begin{tabular}{ccccrrrr}
\toprule
$p$             & $m$ & $n$                            & \multicolumn{2}{c}{time} & memory \\
                                                         \cmidrule(r){4-5}
                &     &                                & (1 core)         & (8 cores)         &             \\
\midrule
\multicolumn{6}{l}{Rader's algorithm with one prime} \\
2147483579 & 1  & 1073741789                           &      400.5       &     56.6          &   36.1      \\
2147477873 & 8  & 134217367                            &      297.4       &     38.9          &    8.0      \\
\midrule
\multicolumn{6}{l}{Rader's algorithm with two primes} \\
2147483543 & 1  & $3137 \cdot 342283$                  &      435.9       &     58.9          &   36.1      \\
2147482577 & 8  & $953 \cdot 140837$                   &      324.4       &     42.2          &    8.1      \\
\midrule
\multicolumn{6}{l}{Rader's algorithm with three primes} \\
2147482367 & 1  & $281 \cdot 1319 \cdot 2897$          &      451.1       &     61.1          &   36.2      \\
2147466449 & 8  & $61 \cdot 337 \cdot 6529$            &      404.9       &     52.5          &    8.1      \\
\midrule
\multicolumn{6}{l}{Rader's algorithm with four primes} \\
2147483399 & 1  & $19 \cdot 31 \cdot 1019 \cdot 1789$  &      492.0       &     64.4          &   37.8      \\
2147453873 & 8  & $17 \cdot 19 \cdot 193 \cdot 2153$   &      417.8       &     54.0          &    8.4      \\
\midrule
\multicolumn{6}{l}{Umbrella algorithm} \\
2147478659 & 1  & $17^2 \cdot 107 \cdot 2671$          &      741.8       &     97.4          &   60.1      \\
2147470673 & 8  & $59^2 \cdot 38557$                   &      500.8       &     66.3          &   11.1      \\
\bottomrule
\end{tabular}
\label{tab:performance}
\end{table}

\section{Kummer--Vandiver conjecture and cyclotomic invariants}
\label{sec:checks}

For verifying the Kummer--Vandiver conjecture, we used essentially the same method as in \cite{BH-irregular} and several previous papers.
This requires computing a certain quantity $V_{p,r}$ (see \cite[\S4]{BH-irregular}) for each irregular pair $(p, r)$, and checking that it does not satisfy a certain congruence.
(Note that there is an error in the formula for $V_{p,r}$ in \cite{BH-irregular} and also in \cite{BCEMS-twelve}; the exponent $c$ should be $c/2$.
Apparently, this error was not propagated to the source code of the programs described in those papers.)

To check the Iwasawa invariants for $p$, it suffices to examine a few congruences for each associated irregular index $r$.
Let
 \[ s_{p,r} = \frac1p \sum_{a=1}^{(p-1)/2} a^{r-1} \in \ZZ, \qquad t_{p,r} = \frac1p \sum_{a=1}^{(p-1)/2} a^{p+r-2} \in \ZZ. \]
If for all irregular indices $r$ the incongruences
 \[ 2^r \not\equiv 1, \qquad s_{p,r} \not\equiv 0, \qquad s_{p,r} \not\equiv t_{p,r}, \qquad (2-r) s_{p,r} \not\equiv (1-r) t_{p,r} \]
all hold modulo $p$, then $\lambda_p = \nu_p = i_p$.
The correctness of this test is an immediate consequence of Criterion 1 and Criterion 3 of \cite[\S3]{EM-cyclotomic}.
(Note that this test is stated incorrectly in \cite{BH-irregular}: the second condition is missing and there is an error in the last condition.)
For $p < 2^{31}$, these incongruences were found to be satisfied for all but five irregular pairs.
The exceptional pairs were
\begin{multline*}
(130811, 52324), (599479, 359568), (2010401, 1234960), \\
 (355011619, 280274852), (358350581, 232032460).
\end{multline*}
For these pairs, the first condition $2^r \not\equiv 1 \pmod p$ failed, and we applied instead the alternative test furnished by Criterion 2 and Criterion~4 of \cite{EM-cyclotomic}.

For each irregular pair, these tests can be carried out in $O(p)$ word operations, using the methods described in \cite{BH-irregular} and \cite{BCEMS-twelve}.
We carried this out using two separate C implementations.
The values of $V_{p,r}$, $s_{p,r}$ and $t_{p,r}$ computed in each run were identical, and are available on request.
One implementation was written by the second author, based on code used in \cite{BH-irregular}, updated to work for the larger primes.
The second implementation was written by the first author from scratch, working independently and sharing no code.

\section{Hardware}
\label{sec:hardware}

The main irregular prime computation was performed over a period of about ten months, starting in late 2012.
We used a number of systems, mainly traditional compute clusters, whose characteristics are summarised in Table \ref{tab:clusters}.
Some of the systems changed configuration during our computations, so for some of the table entries we give representative values.
The main part of the code was written entirely in C, except for a small amount of inline assembly code for the modular arithmetic.
We used various versions of the GCC compiler on different systems.

\textsc{Katana} is the main scientific computation system of the Faculty of Science at the University of New South Wales (UNSW).
Its predecessor \textsc{Tensor} was retired in 2014.
\textsc{Condor} is not really a compute cluster \emph{per se}; that row in the table refers to jobs that were run using idle cycles (especially overnight) on desktop machines in the School of Mathematics at UNSW, managed by a Condor server \cite{condor}.
\textsc{Bowery}, \textsc{Union Square} and \textsc{Cardiac} were the primary HPC systems at New York University (NYU), until they were retired in 2014.
\textsc{Vayu} was the peak system of the National Computational Infrastructure (NCI) facility, funded by the Australian Government and hosted at the Australian National University (ANU), until its retirement in September 2013.
It was replaced by \textsc{Raijin} in June 2013.
Finally, \textsc{Orange} is a system run by Intersect New South Wales; it went into production in early 2013, and the second author was fortunate to obtain access as one of the early test users.
All of these systems are built from contemporary scientific-grade Intel processors, with the exception of \textsc{Cardiac} which used AMD chips, and the \textsc{Condor} machines which contain generally cheaper versions of the Intel hardware.

\begin{table}[h]
\caption{System characteristics}
\begin{tabular}{llrrrr}
\toprule
Location      & Cluster             & Core-hours      & Total  & Cores per & RAM per   \\
              &                     & ($\times 1000$) & cores  & node      & node (GB) \\
\midrule
UNSW          & \textsc{Katana}     &   2004          &  1280  & 12         & 24--144  \\
              & \textsc{Tensor}     &    819          &   336  & 8          & 16--24   \\
              & \textsc{Condor}     &   1840          &   636  & 4          & 16       \\
\midrule
NYU           & \textsc{Bowery}     &   1085          &  2528  & 12         & 24--96   \\
              & \textsc{Union Sq.}  &    383          &   584  & 8          & 16--32   \\
              & \textsc{Cardiac}    &    539          &  1264  & 16         & 32       \\
\midrule
NCI           & \textsc{Vayu}       &     94          & 11936  & 8          & 24       \\
              & \textsc{Raijin}     &    496          & 57472  & 16         & 32--128   \\
\midrule
Intersect     & \textsc{Orange}     &   1370          &  1660  & 16         & 64--256  \\
\bottomrule
\end{tabular}
\label{tab:clusters}
\end{table}

% more detailed breakdown:
%            irregular    correctness    vandiver 1   vandiver 2    cyclotomic 1   cyclotomic 2
%Katana           1383                        148          129            104           240
%Tensor            662                         46           57             54
%Condor           1840
%Bowery            777                         69           69             64           106
%USQ               243                         38           46             56
%Cardiac           424                         34           45             36
%Vayu               94
%Raijin            281             15          40                          59           101
%Orange           1370

As mentioned in Section \ref{sec:irregular}, the memory usage varies considerably depending on which algorithm is deployed.
For example, for $p$ near $2^{31}$, we might need as little as 5 GB or as much as 60 GB.
We grouped our jobs mainly by memory requirements, and allocated them to the various systems depending partly on the available RAM per node and partly on system load.
The systems are shared among many users, and jobs must be submitted using a batch system.
This process was partially automated by a collection of simple Python scripts, but also relied on some manual oversight.

The ``core hours'' column shows the total contribution of each system, including the irregular prime computations, both runs of the Kummer--Vandiver and cyclotomic invariant checks, and a further correctness test (discussed below).
For example, a node with two quad-core CPUs would be listed as contributing 8 core hours for each hour of wall time.
Because the different systems contain different hardware, these numbers are not directly comparable between the rows of the table, but they do give a rough idea of scale.

\section{Correctness}
\label{sec:correctness}

Any computation is susceptible to errors; in a computation of this magnitude it would be a great surprise if nothing went wrong.
Consequently, we took careful precautions, similar to those deployed in \cite{BH-irregular}, to maximise the chance of detecting any problems.

As noted above, the Kummer--Vandiver and cyclotomic checks were executed twice, using independently developed code.
However, for most primes, the expensive determination of the irregular indices was performed only once.
During this step we computed the checksum $C_p = \sum_{r=0}^{p-3} 2^r (r+1) B_r \pmod p$,
and we verified that $C_p \equiv -4 \pmod p$ (see \cite{BCS-million}).
Notice that this test depends on every single $B_r \pmod p$ that we have computed.
Suppose that one of the $B_r \pmod p$ is computed incorrectly, and that we model this by assuming that $C_p$ takes on a random value modulo $p$; then there is still only a $1/p$ chance that we fail to notice the error.
Even if we committed an error for \emph{every} prime between $163\,577\,856$ (the previous search bound in \cite{BH-irregular}) and $2^{31} = 2\,147\,483\,648$, the number of errors that we expect to fail to detect is only about $\log \log(2^{31}) - \log \log(163\,577\,856) \approx 0.128$.
For this reason, we are reasonably confident that we did detect every error that occurred.

Indeed, a number of errors \emph{were} detected.
The consumer-grade machines in the \textsc{Condor} pool tended to have lower quality RAM, and on a handful of them the checksum test would reliably fail several times a day.
The other systems had high-quality error-correcting RAM modules, and we did not detect any errors on them except for one problematic node on \textsc{Katana}.
If any machine exhibited even a single checksum error, we excluded it from all computations and reprocessed all primes that had been handled on that machine.

A weakness of the above checksum test is that it cannot be checked after the fact, because we discard most of the computed residues.
As a backup, for each prime $p$ we store ten pairs $(r, B_r \bmod p)$, sorted by the value of $B_r \pmod p$.
In other words, we first store those $r$ for which $B_r = 0 \pmod p$ (the irregular indices), followed by those $r$ for which $B_r = 1 \pmod p$, and so on, until we have ten values.
The total compressed size of this auxiliary data set is 5.1 GB.
Later, we checked each of these saved values, using an algorithm that directly computes $B_r \pmod p$ for a single $r$ in $O(p)$ word operations \cite{Har-bernmm}.
This computation ran for about 15{,}000 hours on \textsc{Raijin} and encountered no discrepancies.

\section{Historical remarks}
\label{sec:history}

Prior to 1991, all published searches for irregular primes had been based on quadratic-time algorithms, that is, algorithms that require $p^{2+o(1)}$ bit operations to determine the irregular indices for $p$.
The last paper to work under this regime was \cite{TW-bernoulli} in 1987, in which the irregular primes up to $150\,000$ were found.

At this point there was a switch to quasi-linear algorithms, which perform the same task in $p^{1+o(1)}$ bit operations.
In 1991, Sompolski described the ``power series method'' in his Ph.D. thesis \cite{Som-fermat}, which uses FFTs and Newton's method to compute the first $p + O(1)$ terms of the power series $x/(e^x - 1) = \sum_{r=0}^\infty B_r x^r / r!$ modulo $p$.
Buhler and Crandall were working on almost identical techniques simultaneously; they and Sompolski became aware of each other's work, and together published a search going up to $10^6$ \cite{BCS-million}.

The first published quasi-linear algorithm actually appeared slightly earlier, in 1988.
Chellali's algorithm \cite{Che-bernoulli} was based on the inversion of the same power series, and explicitly uses FFTs, but it was asymptotically inferior to the algorithm sketched in the previous paragraph by a factor of $\log p$, because of the use of a less efficient method for power series inversion.
It does not appear that Chellali attempted to use this algorithm in any large-scale computations.

However, there is considerable evidence that A.~O.~L.~Atkin (Sompolski's Ph.D. advisor) had already discovered at least one, and possibly several, quasi-linear algorithms about a decade earlier.
We find only one brief mention of this in print, at the end of section II.B.1 of Sompolski's thesis.
Sompolski has confirmed (personal communication, 2015) that Atkin apparently knew of the ``power series method'' earlier in the 1980s, and possibly even the late 1970s.
Sompolski also told us that Atkin advocated the use of FFTs for these series manipulations, and that he had the impression that Atkin had actually tried this and found it to be efficient.
Certainly the series inversion algorithm based on Newton's method and FFTs was well known by then (see for example \cite{Kun-reciprocals}).

We also discussed this question with Samuel Wagstaff, Jr. (personal communication, 2016).
According to Wagstaff, Atkin already had an implementation of a quasi-linear algorithm in 1978.
Atkin never explained his algorithm in full, and did not make his source code available, although Wagstaff did have several conversations with him around that time.
Wagstaff's recollection is that Atkin's algorithm probably involved computing the coefficients of the polynomial $(x-1)(x-2) \cdots (x-(p-1))$ modulo $p^2$, and then used the Newton--Girard identities to deduce the power sums $1^k + 2^k + \cdots + (p-1)^k \pmod{p^2}$ for the relevant $k$, which in turn yields information about the Bernoulli numbers modulo $p$.
It is unclear how Atkin accomplished all of this in quasi-linear time, or whether FFTs were involved, although Wagstaff believes that the theoretical running time was $O(p (\log p)^m)$ word operations for some small $m$.
Wagstaff also reports that in 1978 he was able to run Atkin's program for primes in the range $10^6 < p < 10^7$, and observed that it was much faster than Wagstaff's own quadratic-time implementation, and that its performance appeared to be quasi-linear.

\section*{Acknowledgments}

The authors thank UNSW, NYU, NCI and Intersect for the provision of computing resources, and in particular Martin Thompson (UNSW) and Joachim Mai (Intersect) for their technical assistance.
The authors are indebted to Robert Sompolski and Samuel Wagstaff, Jr., for their insights into Atkin's unpublished work on computing irregular primes.
The authors thank Zdravko Botev for useful discussions on statistics, Maike Massierer for help in translating portions of Kummer's work, and Joe Buhler for his comments on a draft manuscript.
The first author was supported by DFG Priority Project SPP 1489 and EPSRC grant EP/G004870/1.
The second author was supported by the Australian Research Council, grants DE120101293 and DP150101689.
The third author was supported by an Australian Postgraduate Award.

\bibliographystyle{amsplain}
\bibliography{twobillion}

\providecommand{\bysame}{\leavevmode\hbox to3em{\hrulefill}\thinspace}
\providecommand{\MR}{\relax\ifhmode\unskip\space\fi MR }
% \MRhref is called by the amsart/book/proc definition of \MR.
\providecommand{\MRhref}[2]{%
  \href{http://www.ams.org/mathscinet-getitem?mr=#1}{#2}
}
\providecommand{\href}[2]{#2}
\begin{thebibliography}{10}

\bibitem{Blu-dft}
L.~Bluestein, \emph{A linear filtering approach to the computation of discrete
  {F}ourier transform}, Audio and Electroacoustics, IEEE Transactions on
  \textbf{18} (1970), no.~4, 451--455.

\bibitem{BCEM-four}
J.~P. Buhler, R.~E. Crandall, R.~Ernvall, and T.~Mets{\"a}nkyl{\"a},
  \emph{Irregular primes and cyclotomic invariants to four million}, Math.
  Comp. \textbf{61} (1993), no.~203, 151--153. \MR{1197511 (93k:11014)}

\bibitem{BCEMS-twelve}
J.~P. Buhler, R.~E. Crandall, R.~Ernvall, T.~Mets{\"a}nkyl{\"a}, and M.~A.
  Shokrollahi, \emph{Irregular primes and cyclotomic invariants to 12 million},
  J. Symbolic Comput. \textbf{31} (2001), no.~1-2, 89--96. \MR{1806208
  (2001m:11220)}

\bibitem{BCS-million}
J.~P. Buhler, R.~E. Crandall, and R.~W. Sompolski, \emph{Irregular primes to
  one million}, Math. Comp. \textbf{59} (1992), no.~200, 717--722. \MR{1134717
  (93a:11106)}

\bibitem{BH-irregular}
J.~P. Buhler and D.~Harvey, \emph{Irregular primes to 163 million}, Math. Comp.
  \textbf{80} (2011), no.~276, 2435--2444. \MR{2813369 (2012j:11243)}

\bibitem{Che-bernoulli}
M.~Chellali, \emph{Acc\'el\'eration de calcul de nombres de {B}ernoulli}, J.
  Number Theory \textbf{28} (1988), no.~3, 347--362. \MR{932380 (89h:05002)}

\bibitem{EM-cyclotomic}
R.~Ernvall and T.~Mets{\"a}nkyl{\"a}, \emph{Cyclotomic invariants for primes
  between {$125\,000$} and {$150\,000$}}, Math. Comp. \textbf{56} (1991),
  no.~194, 851--858. \MR{1068819 (91h:11157)}

\bibitem{Goo-fourier}
I.~J. Good, \emph{The interaction algorithm and practical {F}ourier analysis},
  J. Roy. Statist. Soc. Ser. B \textbf{20} (1958), 361--372. \MR{0102888 (21
  \#1674)}

\bibitem{Har-cachetft}
D.~Harvey, \emph{A cache-friendly truncated {FFT}}, Theoret. Comput. Sci.
  \textbf{410} (2009), no.~27-29, 2649--2658. \MR{2531107 (2010g:68327)}

\bibitem{Har-bernmm}
\bysame, \emph{A multimodular algorithm for computing {B}ernoulli numbers},
  Math. Comp. \textbf{79} (2010), no.~272, 2361--2370. \MR{2684369
  (2011h:11019)}

\bibitem{Har-ntt}
\bysame, \emph{Faster arithmetic for number-theoretic transforms}, J. Symbolic
  Comput. \textbf{60} (2014), 113--119. \MR{3131382}

\bibitem{Joh-vanishing}
W.~Johnson, \emph{On the vanishing of the {I}wasawa invariant {$\mu _{p}$} for
  {$p<8000$}}, Math. Comp. \textbf{27} (1973), 387--396. \MR{0384748 (52
  \#5621)}

\bibitem{Joh-irregular}
\bysame, \emph{Irregular primes and cyclotomic invariants}, Math. Comp.
  \textbf{29} (1975), 113--120, Collection of articles dedicated to Derrick
  Henry Lehmer on the occasion of his seventieth birthday. \MR{0376606 (51
  \#12781)}

\bibitem{Kob-fermat}
V.~V. Kobelev, \emph{A proof of {F}ermat's theorem for all prime exponents less
  than {$5500$}.}, Dokl. Akad. Nauk SSSR \textbf{190} (1970), 767--768.
  \MR{0258717}

\bibitem{Kum-allgemeiner}
E.~E. Kummer, \emph{Allgemeiner {B}eweis des {F}ermatschen {S}atzes, da\ss\ die
  {G}leichung {$x^\lambda+y^\lambda=z^\lambda$} durch ganze {Z}ahlen unl\"osbar
  ist, f\"ur alle diejenigen {P}otenz-{E}xponenten {$\lambda$} welche ungerade
  {P}rimzahlen sind und in den {Z}\"ahlern der ersten {$\frac12(\lambda-3)$}
  {B}ernoullischen zahlen als {F}actoren nicht vorkommen}, J. Reine Angew.
  Math. \textbf{40} (1850), 130--138. \MR{1578681}

\bibitem{Kum-memoire}
\bysame, \emph{M\'emoire sur la th\'eorie des nombres complexes compos\'es de
  racines de l'unit\'e et de nombres entiers}, J. Math. Pure Appl \textbf{16}
  (1851), 377--498.

\bibitem{Kum-1874}
\bysame, \emph{{Ueber diejenigen Primzahlen $\lambda$, f\"ur welche die
  Klassenzahl der aus $\lambda^{\text{ten}}$ Einheitswurzeln gebildeten
  complexen Zahlen durch $\lambda$ theilbar ist.}}, {Berl. Monatsber.} (1874),
  239--248 (German).

\bibitem{Kun-reciprocals}
H.~T. Kung, \emph{On computing reciprocals of power series}, Numer. Math.
  \textbf{22} (1974), 341--348. \MR{0351045}

\bibitem{Lan-cyclotomic-combined}
S.~Lang, \emph{Cyclotomic fields {I} and {II}}, second ed., Graduate Texts in
  Mathematics, vol. 121, Springer-Verlag, New York, 1990, With an appendix by
  Karl Rubin. \MR{1029028 (91c:11001)}

\bibitem{Leh-automation}
D.~H. Lehmer, \emph{Automation and pure mathematics}, Applications of Digital
  Computers (W.~F. Freiberger and W.~Prager, eds.), Ginn, Boston, Mass, 1963.

\bibitem{LLV-FLT}
D.~H. Lehmer, E.~Lehmer, and H.~S. Vandiver, \emph{An application of high-speed
  computing to {F}ermat's last theorem}, Proc. Nat. Acad. Sci. U. S. A.
  \textbf{40} (1954), 25--33. \MR{0061128 (15,778f)}

\bibitem{Rad-prime}
C.~M. Rader, \emph{Discrete {Fourier} transforms when the number of data
  samples is prime}, Proc. IEEE \textbf{56} (1968), no.~6, 1107--1108.

\bibitem{SNV-FLT}
J.~L. Selfridge, C.~A. Nicol, and H.~S. Vandiver, \emph{Proof of {F}ermat's
  last theorem for all prime exponents less than {$4002$}}, Proc. Nat. Acad.
  Sci. U.S.A. \textbf{41} (1955), 970--973. \MR{0072892}

\bibitem{Sho-eight}
M.~A. Shokrollahi, \emph{Computation of irregular primes up to eight million
  (preliminary report)}, Tech. Report TR-96-002, International Computer Science
  Institute, Berkeley, 1996.

\bibitem{Som-fermat}
R.~W. Sompolski, \emph{The second case of {F}ermat's last theorem for fixed
  irregular prime exponents}, Ph.D. thesis, University of Illinois at Chicago,
  1991.

\bibitem{SV-irregular}
E.~T. Stafford and H.~S. Vandiver, \emph{Determination of some properly
  irregular cyclotomic fields}, Proceedings of the National Academy of Sciences
  \textbf{16} (1930), no.~2, 139--150.

\bibitem{TW-bernoulli}
J.~W. Tanner and S.~Wagstaff, Jr., \emph{New congruences for the {B}ernoulli
  numbers}, Math. Comp. \textbf{48} (1987), no.~177, 341--350. \MR{866120
  (87m:11017)}

\bibitem{condor}
The~Condor team, \emph{Condor software package},
  http://research.cs.wisc.edu/htcondor/.

\bibitem{Tho-physics}
L.~H. Thomas, \emph{Using a computer to solve problems in physics},
  Applications of Digital Computers (1963), 44--45.

\bibitem{vdH-TFT-apps}
J.~van~der Hoeven, \emph{The truncated {F}ourier transform and applications},
  I{SSAC} 2004, ACM, New York, 2004, pp.~290--296. \MR{MR2126956}

\bibitem{Van-FLT}
H.~S. Vandiver, \emph{On {B}ernoulli's numbers and {F}ermat's last theorem},
  Duke Math. J. \textbf{3} (1937), no.~4, 569--584. \MR{1546011}

\bibitem{Van-attack}
\bysame, \emph{Examination of methods of attack on the second case of
  {F}ermat's last theorem}, Proc. Nat. Acad. Sci. U. S. A. \textbf{40} (1954),
  732--735. \MR{0062758}

\bibitem{Wag-irregular}
S.~Wagstaff, Jr., \emph{The irregular primes to {$125000$}}, Math. Comp.
  \textbf{32} (1978), no.~142, 583--591. \MR{0491465 (58 \#10711)}

\bibitem{Was-cyclotomic}
L.~C. Washington, \emph{Introduction to cyclotomic fields}, second ed.,
  Graduate Texts in Mathematics, vol.~83, Springer-Verlag, New York, 1997.
  \MR{1421575 (97h:11130)}

\end{thebibliography}

\end{document}